\pdfoutput=1
\RequirePackage{ifpdf}
\ifpdf 
\documentclass[pdftex]{sigma}
\else
\documentclass{sigma}
\fi

\numberwithin{equation}{section}

\newtheorem{Theorem}{Theorem}[section]
\newtheorem*{Theorem*}{Theorem}
\newtheorem{Corollary}[Theorem]{Corollary}

\newtheorem{Proposition}[Theorem]{Proposition}
 { \theoremstyle{definition}

\newtheorem{Remark}[Theorem]{Remark} }

\DeclareMathOperator{\End}{End}
\DeclareMathOperator{\Int}{Int}
\DeclareMathOperator{\Ker}{Ker}
\DeclareMathOperator{\vol}{vol}
\DeclareMathOperator{\Tr}{Tr}
\DeclareMathOperator{\Id}{Id}
\DeclareMathOperator{\prim}{prim}

\newcommand{\ii}{\sqrt{-1}}
\newcommand{\adjoint}{\dagger}
\newcommand{\deRham}{\mathrm{dR}}
\newcommand{\Rumin}{\mathrm{R}}
\newcommand{\RN}{\cE}

\newcommand{\bR}{\mathbb{R}}
\newcommand{\bC}{\mathbb{C}}
\newcommand{\cE}{\mathcal{E}}
\newcommand{\rest}[1]{\big\rvert_{#1}}

\begin{document}
\allowdisplaybreaks

\newcommand{\arXivNumber}{2009.03276}

\renewcommand{\PaperNumber}{091}

\FirstPageHeading

\ShortArticleName{Ray--Singer Torsion and the Rumin Laplacian on Lens Spaces}

\ArticleName{Ray--Singer Torsion and the Rumin Laplacian\\ on Lens Spaces}

\Author{Akira KITAOKA}

\AuthorNameForHeading{A.~Kitaoka}

\Address{Graduate School of Mathematical Sciences, The University of Tokyo, Japan}
\Email{\href{akira5kitaoka@gmail.com}{akira5kitaoka@gmail.com}}
\URLaddress{\url{https://akira5kitaoka.github.io/Akira5Kitaoka-en.github.io/}}

\ArticleDates{Received May 19, 2022, in final form November 14, 2022; Published online November 28, 2022}

\Abstract{We express explicitly the analytic torsion functions associated with the Rumin complex on lens spaces in terms of the Hurwitz zeta function. In particular, we find that the functions vanish at the origin and determine the analytic torsions. Moreover, we have a~formula between this torsion and the Ray--Singer torsion.}

\Keywords{analytic torsion; Rumin complex; CR geometry; contact geometry}

\Classification{58J52; 32V20; 53D10; 43A85}

\section{Introduction}

Let $(M,H)$ be a compact contact manifold of dimension $2n+1$
and $E$ be the flat vector bundle with unitary holonomy on $M$.
Rumin~\cite{Rumin-94} introduced a complex
$(\cE ^\bullet (M ,E) , {\rm d}_{\Rumin }^\bullet )$,
which is
a~subquotient of the de Rham complex of $E$.
A specific feature of the complex is that
the operator $D = {\rm d}_{\Rumin }^n \colon \cE ^{n} (M,E) \to \cE ^{n+1} (M,E) $ in `middle degree' is a second-order,
while
 ${\rm d}_{\Rumin }^k \colon \cE ^{k} (M,E) \to \cE ^{k+1} (M,E) $ for $k \not = n$ are first order which are induced by
the exterior derivatives.
Let $a_k = 1 / \sqrt{|n-k|}$ for $k \not = n$ and $a_n = 1$.
Then, $(\cE ^\bullet (M,E) , {\rm d}_\RN ^\bullet)$,
where ${\rm d}_\RN ^k = a_k {\rm d}_{\Rumin}^k$,
is also a complex.
We call $(\cE ^\bullet (M,E) , {\rm d}_\RN ^\bullet)$ the {\em Rumin complex}.
In virtue of the rescaling, ${\rm d}_\RN ^\bullet$~satisfies K\"{a}hler-type identities on Sasakian manifolds~\cite[equation~(34)]{Rumin-00}, which include the case of lens spaces.

Let $\theta$ be a contact form of $H$
and
$J$ be an almost complex structure on $H$.
Then we may define a Riemann metric $g_{\theta , J}$ on $TM$ by extending the Levi metric ${\rm d} \theta (- ,J-)$ on $H$
(see Section~\ref{sec: The Rumin complex on contact manifolds}).
Following~\cite{Rumin-94}, we define the Rumin Laplacians $\Delta_\RN $ associated with $(\cE ^\bullet (M,E) , {\rm d}_\RN ^\bullet)$ and the metric~$g_{\theta , J}$ by
\[
 \Delta_\RN ^k
 :=
 \begin{cases}
 ( {\rm d}_\RN {\rm d}_\RN {}^{\adjoint })^2
 + ({\rm d}_\RN {}^{\adjoint } {\rm d}_\RN )^2 , & k \not = n, \, n+1 , \\
 ( {\rm d}_\RN {\rm d}_\RN {}^{\adjoint })^2
 + D^{\adjoint } D , & k=n, \\
 D D^{\adjoint }
 + ({\rm d}_\RN {}^{\adjoint } {\rm d}_\RN )^2 , & k =n+1.
 \end{cases}
\]
Rumin showed that $\Delta_\RN $ has discrete eigenvalues with finite multiplicities.

We next introduce the analytic torsion and metric of the Rumin complex $(\cE ^\bullet (M,E) , {\rm d}_\RN ^\bullet)$ by~following~\cite{Bismut-Zhang-92, Kitaoka-19, Rumin-Seshadri-12}.
We define the {\em contact analytic torsion function} associated with $(\cE ^\bullet (M,E),\allowbreak {\rm d}_\RN ^\bullet)$ by
\begin{equation}
 \label{kappa}
 \kappa_\RN ( M, E, g_{\theta , J}) (s)
 := \sum_{k=0}^{n} (-1)^{k+1} (n+1-k) \zeta \big(\Delta_\RN ^k \big) (s)
 ,
\end{equation}
where $\zeta \big(\Delta_\RN ^k \big) (s)$ is the spectral zeta function of $\Delta_\RN ^k$.
Since the Rumin Laplacians $\Delta_\RN $ satisfies the Rockland condition by~\cite[p.~300]{Rumin-94},
the spectral zeta function $\zeta \big(\Delta_\RN ^k \big) (s)$ extends to a meromorphic function
on $\mathbb{C}$ which is holomorphic at zero by~\cite[Section~4]{Ponge-07}.
Here, we use $0^s = 1$.
We define the {\em contact analytic torsion} $ T_\RN $ by
\[
 2 \log T_\RN ( M , E, g_{\theta , J})
 =
 \kappa_\RN ( M , E, g_{\theta , J}) ' (0) .
\]
Let $H^\bullet (\cE ^\bullet , {\rm d}_\RN ^\bullet)$ be the cohomology of the Rumin complex.
We define the contact metric on $\det H^\bullet (\cE ^\bullet , {\rm d}_\RN ^\bullet)$ by
\[
 \| \ \|_\RN (M, E, g_{\theta , J} )
 =
 T_\RN^{-1} ( M , E, g_{\theta , J})
 | \ |_{L^2 (\cE ^\bullet) },
\]
where
the metric
$| \ |_{L^2 (\cE ^\bullet) }$ is induced by $L^2$ metric on $\cE ^\bullet (M, E)$
via identification of the cohomology classes by the harmonic forms on $\cE ^\bullet (M, E)$.

Rumin and Seshadri~\cite{Rumin-Seshadri-12} defined another analytic torsion function $\kappa_{\Rumin}$ from
$(\cE ^\bullet (M,E) , {\rm d}_{\Rumin }^\bullet )$, which is different from $\kappa_\RN $ except in dimension 3.
In dimension 3, they showed that
$\kappa_{\Rumin} (M, E, \allowbreak g_{\theta , J} ) (0)$ is a contact invariant,
that is, independent of the metric $g_{\theta , J}$.
Moreover, on 3-dimensional Sasakian manifolds with $S^1$-action, $\kappa_{\Rumin} (M, E, g_{\theta , J} ) (0) = 0 $.
Here,
Sasakian manifolds with $S^1$-action means Sasakian manifolds whose Reeb vector filed generates the circle action $S^1$.
3-di\-men\-sional Sasakian manifolds with $S^1$-action
are CR Seifert manifolds.
Furthermore, they showed that this analytic torsion and the Ray--Singer torsion $T_{\deRham } (M, E, g_{\theta , J} )$ equal for flat bundles with unitary holonomy on 3-dimensional Sasakian manifolds with $S^1$ action.

To extend the coincidence,
with ${\rm d}_{\RN}$ instead of ${\rm d}_{\Rumin }$,
the author~\cite{Kitaoka-19} showed that $T_{\RN } \big(S^{2n+1}\!, \underline{\bC}, g_{\theta , J} \big) \\= n! T_{\deRham } \big(S^{2n+1}, \underline{\bC}, g_{\theta , J} \big)$ on the standard CR spheres $S^{2n+1} \big({\subset} \bC^{n+1} \big) $,
where $\underline{\bC}$ is the trivial line bundle.
Here the {\em standard CR sphere} is triple $\big(S^{2n+1}, \theta, J\big)$,
where
$\theta$ is given the contact form by $\theta = \sqrt{-1} \big(\bar{\partial} - \partial\big) |z|^2$
and
$J$ is an almost complex structure $J$ induced from the complex structure of $ \bC^{n+1}$.
It is simply denoted by $S^{2n+1}$.
Moreover, Albin and Quan~\cite[Corollary~3 and equation~(4)]{Albin-Quan-20} showed the difference between the Ray--Singer torsion and the contact analytic torsion is given by some integrals of universal polynomials in the local invariants of the metric
on contact manifolds.

In this paper,
we extend this coincidence on lens spaces
and determine explicitly the analytic torsion functions associated with the Rumin complex
in terms of the Hurwitz zeta function.
Let $g_{\mathrm{std}}$ be the standard metric on $S^{2n+1}$
and we note that $g_{\theta , J} = 4 g_{\mathrm{std}}$.
Let $\mu , \nu_1 , \ldots , \nu_{n+1}$ be integers such that the $\nu_j$ are coprime to $\mu$.
Let $\Gamma$ be the subgroup of $\big(S^1\big)^{n+1}$ generated by
\begin{equation}
 \gamma
 = ( \gamma_1 , \ldots , \gamma_{n+1} )
 := \big( \exp \big( 2 \pi \ii \nu_1 / \mu \big), \dots , \exp \big( 2 \pi \ii \nu_{n+1} / \mu \big) \big) .
 \notag
\end{equation}
We denote the lens space by
\[
 K
 :=
 S^{2n+1} / \Gamma
 .
\]

Let $\underline{\bC}$ be the trivial line bundle on $K$.
Fix $u \in \mathbb{Z}$ and consider the unitary representation $\alpha_u \colon \pi_1 (K) = \Gamma \to \mathrm{U} (1)$, defined by
\[
 \alpha_u \big( \gamma^{\ell} \big) := \exp \big( 2 \pi \ii u \ell / \mu \big)
 \qquad \text{for} \ \ell \in \mathbb{Z}
 .
\]
Let $E_{\alpha}$ be the flat vector bundle associated with the unitary representation ${\alpha \colon \pi_1 (K) = \Gamma \! \to \! \mathrm{U} (r)}$, and $E_{\alpha_u} = E_{u}$.
The sections of this bundle correspond to $\alpha_u$-equivariant functions on~$S^{2n+1}$.

Our main result is
\begin{Theorem}
\label{theo:compute-precisely-kappa-on-the-lens-spaces}

Let $K$ be the lens space with
the contact form and the almost complex structure which are induced by the action $\Gamma$ on the standard CR sphere $S^{2n+1}$.
\begin{enumerate}\itemsep=0pt
\item[$(1)$]
 The contact analytic torsion function
 of $(K, \underline{\bC})$
 is given by
 \begin{equation}
 \kappa_\RN ( K , \underline{\bC}, g_{\theta , J} ) (s)
 =
 -(n+1)
 \bigr(
 1
 +
 2^{2s+1}
 \mu^{-2s}
 \zeta (2s)
 \bigl)
 ,
 \label{eq:kappa-on-lens-spaces}
 \end{equation}
 where $\zeta$ is the Riemann zeta function.
 In particular, we have
 \begin{align}
 & \kappa_\RN ( K , \underline{\bC} , g_{\theta , J} ) (0)
 =
 0
 ,
 \label{eq:kappa-at-the-origin-of-the-triv-bundle-on-lens-spaces}
 \\
 & T_\RN ( K , \underline{\bC}, g_{\theta , J} )
 =
 \left(\frac{4 \pi}{\mu}\right)^{n+1}
 .
 \label{eq:Rumin-analytic-torsion-of-the-triv-bundle-on-lens-spaces}
 \end{align}

\item[$(2)$]
 The contact analytic torsion function
 of $(K, E_{u})$
 for $u \in \{ 1, \ldots , \mu -1 \} $
 is given by
 \begin{equation}
 \label{eq:kappa-of-flat-bundles-on-lens-spaces}
 \kappa_\RN ( K , E_{u}, g_{\theta , J} ) (s)
 =
 -2^{2s} \mu^{-2s}
 \sum_{j=1}^{n+1}
 \bigl(
 \zeta ( 2s , A_{\mu} (u {\tau}_j) / \mu )
 +
 \zeta ( 2s , A_{\mu} (- u {\tau}_j) / \mu )
 \bigr)
 ,
 \end{equation}
 where $\zeta (s , a) := \sum_{q=0}^{\infty} (q +a )^{-s}$ is the Hurwitz zeta function for $0 < a \leq 1$, $A_{\mu} ( w )$ is the integer between $1$ and $\mu$
 such that $A_{\mu} ( w ) \equiv w \mod \mu$
 and ${\tau}_j \nu_j \equiv 1 \mod \mu$.
 In particular, we have
 \begin{align}
 & \kappa_\RN ( K , E_{u}, g_{\theta , J} ) (0)
 =
 0
 ,
 \label{eq:kappa-at-the-origin-of-flat-vect-bundle-on-lens-spaces}
 \\
 & T_\RN ( K , E_{u}, g_{\theta , J} )
 =
 \prod_{j=1}^{n+1} \big| {\rm e}^{ 2 \pi \ii u {\tau}_j / \mu } -1 \big|
 .
 \label{eq:Rumin-analytic-torsion-of-flat-vect-bundle-on-lens-spaces}
 \end{align}
\end{enumerate}
\end{Theorem}

The equations \eqref{eq:kappa-on-lens-spaces} and \eqref{eq:kappa-of-flat-bundles-on-lens-spaces} extend the following results of $\kappa_{\RN}$ the spheres to on lens spaces.
Rumin and Seshadri~\cite[Theorem 5.4]{Rumin-Seshadri-12} showed \eqref{eq:kappa-on-lens-spaces} in the case of 3-dimensional lens spaces.
The author~\cite{Kitaoka-19} showed \eqref{eq:kappa-on-lens-spaces} in the case of $\big(S^{2n+1}, \underline{\bC}\big)$ for arbitrary $n$.

From \eqref{eq:kappa-at-the-origin-of-the-triv-bundle-on-lens-spaces} and \eqref{eq:kappa-at-the-origin-of-flat-vect-bundle-on-lens-spaces},
we see that
the metric $\| \ \|_\RN $ on $( K, E_{u} , g_{\theta , J} )$ is invariant under the constant rescaling $\theta \mapsto C \theta$.
The argument is exactly the same as the one in~\cite{Rumin-Seshadri-12}.

In the same way as~\cite{Kitaoka-19},
the fact that the representations determine the eigenvalues of $\Delta_\RN$ cause several cancellations in the linear combination \eqref{kappa}, which significantly simplifies the computation of $\kappa_\RN (s)$.
We cannot get such a simple formula for the contact analytic torsion function $\kappa_{\Rumin}$ of $(\cE ^\bullet (M,E), {\rm d}_{\Rumin}^\bullet)$
for dimensions higher than $3$.

Let us compare the contact analytic torsion with the Ray--Singer torsion on lens spaces.
Ray~\cite{Ray-70} showed that for $u (= 1 ,\ldots , \mu -1)$
\[
 T_{\deRham} ( K , E_{u}, 4 g_{\mathrm{std}} )
 =
 \prod_{j=1}^{n+1} \big| {\rm e}^{ 2 \pi \ii u {\tau}_j / \mu } -1 \big|
 .
\]
Weng and You~\cite{Weng-You-96} calculate the Ray--Singer torsion on spheres.
We extend their results for the trivial bundle on lens spaces:
\begin{Proposition}
 \label{prop:Ray--Singer_torsion_of_the_trivial_bundle_on_lens_spaces}
 In the setting of Theorem~{\rm \ref{theo:compute-precisely-kappa-on-the-lens-spaces}},
 we have
 \[
 T_{\deRham} ( K, \underline{\bC} , 4 g_{\mathrm{std}})
 = \frac{ (4 \pi)^{n+1}}{ n! \mu^{n+1}} .
 \]
\end{Proposition}

The metric $4 g_{\mathrm{std}}$ agrees with the metric $g_{\theta , J}$ defined from the contact form $\theta = \sqrt{-1} \big(\bar{\partial} - \partial\big) |z|^2$.
Since the cohomology of $(\cE ^\bullet (M,E) , {\rm d}_\RN ^\bullet)$
coincides with that of $ (\Omega^\bullet (M,E) , {\rm d} )$ (e.g.,~\cite[p.~286]{Rumin-94}), there is a natural isomorphism
$\det H^\bullet (\cE ^\bullet (M,E) , {\rm d}_\RN ^\bullet)
\cong
\det H^\bullet (\Omega^\bullet (M,E) , {\rm d} )
$, which turns out to be isometric for the $L^2$ metrics.
Therefore \eqref{eq:Rumin-analytic-torsion-of-the-triv-bundle-on-lens-spaces} and
\eqref{eq:Rumin-analytic-torsion-of-flat-vect-bundle-on-lens-spaces} give
\begin{Corollary}
\label{cor:Rumin-and-Ray--Singer-analytic-torsion-on-lens-spaces}
In the setting of Theorem~{\rm \ref{theo:compute-precisely-kappa-on-the-lens-spaces}},
for all unitary holonomy $\alpha \colon \pi_{1}(K) \to \mathrm{U} (r)$,
we have
\begin{align*}
& T_{\RN} ( K, E_{\alpha} , g_{\theta , J})
 =
 n!^{\dim H^0 ( K , E_{\alpha}) }
 T_{\deRham} ( K, E_{\alpha} , g_{\theta , J})
 ,
 \\
& \| \ \|_\RN ( K, E_{\alpha} , g_{\theta , J})
=
 n!^{- \dim H^0 ( K , E_{\alpha}) }
 \| \ \|_{\deRham} ( K, E_{\alpha} , g_{\theta , J})
 ,
\end{align*}
via the isomorphism $\det H^\bullet (\cE ^\bullet (M,E_{\alpha} ) , {\rm d}_\RN ^\bullet)
\cong\det H^\bullet (\Omega^\bullet (M,E_{\alpha}) , {\rm d} )$.
\end{Corollary}

The paper is organized as follows.
In Section~\ref{sec: The Rumin complex on the spheres},
we recall the definition and properties of the Rumin complex on $S^{2n+1}$.
In Section~\ref{sec:on_flat_vector_bundles},
we calculate the contact analytic torsion function $\kappa_\RN $ of flat vector bundles on lens spaces.
In Section~\ref{sec:Ray--Singer_torsion_on_the_trivial_vector_bundle},
we compute the Ray--Singer torsion $T_{\deRham} $ of the trivial vector bundle.
In Section~\ref{sec:Ray--Singer_torsion_and_contact_analytic_torsion},
we compare the Ray--Singer torsion and the contact analytic torsion.

\section{The Rumin complex}
\label{sec: The Rumin complex on the spheres}

\subsection{The Rumin complex on contact manifolds}
\label{sec: The Rumin complex on contact manifolds}

We call $(M,H)$ an orientable contact manifold of dimension $2n+1$
if $H$ is a subbundle of $TM$ of codimension $1$
and
there exists a 1-form $\theta$,
called a contact form,
such that
 ${\Ker (\theta \colon TM \! \to \! \bR ) = H}$
and $\theta \wedge ( {\rm d} \theta )^n \not = 0$.
The Reeb vector field of $\theta$ is the unique vector field $T$ satisfying $\theta (T) =1$ and $\Int_T {\rm d} \theta =0$, where $\Int_T$ is the interior product with respect to~$T$.

For $H$ and $\theta$,
we call $J \in \End (TM)$
an almost complex structure associated with $\theta$ if
$J^2 = - \Id $ on $H$, $JT=0$, and the Levi form ${\rm d} \theta ( - , J - )$ is positive definite on~$H$.
Given $\theta $ and $J$, we define a Riemannian metric $g_{\theta , J}$ on $TM$ by
\begin{align*}
g_{\theta , J} (X,Y) := {\rm d} \theta (X , JY) + \theta (X) \theta (Y) \qquad \text{for} \ X,Y \in TM .
\end{align*}
Let $*$ be the Hodge star operator on $\wedge^\bullet T^*M$ with respect to $g_{\theta , J}$.

Let $M$ be a manifold, $\widetilde M$ be the universal covering of $M$,
$\pi_1 (M)$ be the fundamental group of $M$.
For each unitary representation $\alpha \colon \pi_1 (M) \to \mathrm{U} (r )$,
we denote the flat vector bundle associated with $\alpha$ by
\[
 E_{\alpha} := \widetilde M \times_{\alpha} \bC^r \to M .
\]
Let $\nabla_{\alpha}$ be the flat connection on $E_{\alpha}$ induced from the trivial connection on $\widetilde M \times \bC^r$,
and ${\rm d}^{\nabla_{\alpha}}$ be the exterior covariant derivative of $\nabla_{\alpha}$.

The Rumin complex~\cite{Rumin-94} is defined on contact manifolds as follows.
We set $L := {\rm d} \theta \wedge $
and $\Lambda := *^{-1} L *$,
which is the adjoint operator of $L$ with respect to the metric $g_{\theta , J}$ at each point.
We set
\begin{align*}
& \textstyle \bigwedge_{\prim}^{k} H^*:=
\big\{
v \in \textstyle \bigwedge^{k} H^*
\, \big| \,
\Lambda v = 0
\big\},
\\
& \textstyle \bigwedge_{L}^{k} H^*
:=
\big\{
v \in \textstyle \bigwedge^{k} H^*
\, \big| \,
L v = 0
\big\},
\\
& \cE ^k (M , E_{\alpha})
:=
\begin{cases}
 C^\infty \big( M, \bigwedge_{\prim}^{k} H^* \otimes E_{\alpha} \big) , & k \leq n ,\\
 C^\infty \big( M, \theta \wedge\bigwedge_{L}^{k-1} H^* \otimes E_{\alpha} \big) , & k \geq n+1 .
\end{cases}
\end{align*}
We embed $H^*$ into $T^*M$ as the subbundle $ \{ \phi \in T^* M \,|\, \phi (T) = 0 \}$ so that
we can regard
\begin{align*}
\Omega_H^{k} (M, E_{\alpha}) := C^{\infty} \big( M, \textstyle \bigwedge^{k} H^* \otimes E_{\alpha} \big)
\end{align*}
as a subspace of $\Omega^k (M)$, the space of $k$-forms.
We define ${\rm d}_b \colon \Omega_H^{k} (M, E_{\alpha}) \to \Omega_H^{k+1} (M, E_{\alpha})$ by
\[
{\rm d}_b \phi := {\rm d}^{\nabla_{\alpha}} \phi - \theta \wedge \big(\Int_T {\rm d}^{\nabla_{\alpha}} \phi \big)
,
\]
and then $D \colon \cE ^n (M,E_{\alpha}) \to \cE ^{n+1} (M,E_{\alpha}) $ by
\[
D = \theta \wedge \big( \mathcal{L}_T + {\rm d}_b L^{-1} {\rm d}_b \big),
\]
where $\mathcal{L}_T$ is the Lie derivative with respect to $T$,
and we use the fact that $L \colon \textstyle \bigwedge^{n-1} H^* \to\textstyle \bigwedge^{n+1} H^*$ is an isomorphism.

Let $P \colon \textstyle \bigwedge^{k} H^* \to \textstyle \bigwedge_{\prim}^{k} H^* $ be the fiberwise orthogonal projection with respect to $g_{\theta , J}$, which also defines a projection $P \colon \Omega^k (M,E_{\alpha}) \to \cE ^k (M,E_{\alpha}) $.
We set
\[
{\rm d}_{\Rumin }^k
:=
\begin{cases}
P \circ {\rm d}^{\nabla_{\alpha}} & \text{on } \cE ^k (M,E_{\alpha}), \ k \leq n-1, \\
D & \text{on } \cE ^n (M,E_{\alpha}) , \\
{\rm d}^{\nabla_{\alpha}} & \text{on } \cE ^k (M,E_{\alpha}) , \ k \geq n+1.
\end{cases}
\]
Then $(\cE ^\bullet (M,E_{\alpha}), {\rm d}_{\Rumin }^\bullet)$ is a complex.
Let ${\rm d}_\RN ^{k} = a_k {\rm d}_{\Rumin }^{k}$, where
$a_k = 1 / \sqrt{|n-k|}$ for $k \not = n$ and
$a_n = 1$.
We call $(\cE ^\bullet (M,E_{\alpha}), {\rm d}_\RN ^\bullet)$ the Rumin complex.

We define the $L^2$-inner product on $\Omega^{k} (M,E_{\alpha})$ by
\[
(\phi , \psi ) := \int_M g_{\theta , J} (\phi , \psi )\, {\rm d}\vol_{g_{\theta , J}}
\]
and the $L^2$-norm on $\Omega^k (M,E_{\alpha})$ by $\| \phi \| := \sqrt{(\phi , \phi )}$.
Since the Hodge star operator $*$ induces a bundle isomorphism
from $ \textstyle \bigwedge_{\prim}^{k} H^*$ 
to $\theta \wedge \textstyle \bigwedge_{L}^{2n-k} H^*$,
it also induces a map $\cE ^{k} (M,E_{\alpha} ) \to \cE ^{2n+1-k} (M,E_{\alpha}) $.
We note that
\[
 \cE^k ( M , E_{\alpha})
 =
 \big\{
 \phi \in \cE^k \big( \widetilde M , \bC^r \big)
 \, | \,
 t_* \phi = \alpha (t)^{-1} \phi \text{ for } t \in \pi_1 (M)
 \big\}
 .
\]

Let ${\rm d}_\RN ^{\adjoint }$ and $D^{\adjoint }$ denote the formal adjoint of ${\rm d}_\RN $ and $D$, respectively for the $L^2$-inner product.
We define the fourth-order Laplacians $\Delta_\RN $ on $\cE ^k (M,E_{\alpha} )$ by
\[
\Delta_\RN ^k
 :=
\begin{cases}
\big( {\rm d}_\RN ^{k-1} {\rm d}_\RN ^{k-1} {}^{\adjoint }\big)^2 + \big({\rm d}_\RN ^{k} {}^{\adjoint } {\rm d}_\RN ^{k} \big)^2 , & k \not = n , n+1 , \\
\vspace{1mm}
\big( {\rm d}_\RN ^{n-1} {\rm d}_\RN ^{n-1} {}^{\adjoint }\big)^2 + D^{\adjoint } D , & k=n, \\
D D^{\adjoint } + \big({\rm d}_\RN ^{n+1} {}^{\adjoint } {\rm d}_\RN ^{n+1} \big)^2 , & k =n+1.
\end{cases}
\]
We call it the Rumin Laplacian~\cite{Rumin-94}.
Since $*$ and $\Delta_\RN $ commute,
to determine the eigenvalue on $\cE ^\bullet (M,E_{\alpha} ) $,
it is enough to compute them on $\cE ^{k } (M,E_{\alpha})$ for $k \leq n$.

\subsection{The Rumin complex on the CR spheres}

Let $S :=\{ z \in \bC^{n+1} \,|\, |z|^2 =1 \}$
and
$\theta:= \sqrt{-1} \big(\bar{\partial } - \partial \big)|z|^2$.
(We will omit the dimension from~$S^{2n+1}$ for the simplicity of the notation.)
Let $g_{\mathrm{std}}$ be the standard metric on~$S$.
Then, $g_{\theta , J}$ coincides with~$4 g_{\mathrm{std}}$.
With respect to the standard almost complex structure $J$, we decompose the bundles defined in the previous subsection as follows:
\begin{align*}
& H^{* 1,0}:= \big\{ v \in \bC H^* \,|\, Jv = \sqrt{-1} v \big\} , \\
& H^{* 0,1} := \big\{ v \in \bC H^* \,|\, Jv = - \sqrt{-1} v \big\}, \\
& \textstyle \bigwedge^{i, j} H^* := \textstyle \bigwedge^{i} H^{* 1,0} \otimes \textstyle \bigwedge^{j} H^{* 0,1} , \\
& \textstyle \bigwedge_{\prim}^{i, j} H^* := \big\{ \phi \in \textstyle \bigwedge^{i, j} H^* \, \big| \, \Lambda \phi =0 \big\} , \\
%
%
%
%
& \cE ^{i, j} := C^{\infty} \big( S, \textstyle \bigwedge_{\prim}^{i, j} H^* \big).
\end{align*}
We decompose $\cE ^{i,j} $ into a direct sum of irreducible representations of the unitary group~$\mathrm{U} (n+1)$.
Recall that
irreducible representations of $\mathrm{U} (m)$ are parametrized by the highest weight $\lambda = (\lambda_1 , \ldots , \lambda_m ) \in \mathbb{Z}^m$ with $\lambda_1 \geq \lambda_2 \geq \cdots \geq \lambda_m$; the representation corresponding to~$\lambda$ will be denoted by
$ V ( \lambda ) $.
To simplify the notation, we introduce the following notation:
for $a_1, \ldots $, $a_l \in \mathbb{Z}$ and $k_1, \ldots , k_l \in \mathbb{Z}$, $(\underline{a_1}_{k_1}, \dots , \underline{a_l}_{k_l})$ denotes the $k_1 + \cdots + k_l $-tuple
whose first $k_1$ entries are $a_1$, whose next $k_2$ entries are $a_2$, etc.
For example,{\samepage
\[
(\underline{1}_{3}, \underline{0}_{2}, \underline{-1}_{2} ) = (1, 1, 1, 0, 0, -1, -1 ).
\]
We note that $\underline{a}_{1}$ is $a$, and $\underline{a}_{0}$ is the zero tuple.}

In~\cite{Julg-Kasparov-95}, it is shown that
the multiplicity of $V (q, \underline{1}_{j}, \underline{0}_{n-1-i-j}, \underline{-1}_{i}, -p )$ in $\cE ^{s,t}$
is at most one.
Thus we may set
\begin{equation}
\Psi_{(q , j , i , p )}^{(s, t)}
:= \cE ^{s,t}
\cap V (q, \underline{1}_{j}, \underline{0}_{n-1-i-j}, \underline{-1}_{i}, -p )
.
\notag
\end{equation}

\begin{Proposition}[{\cite[Section 4(b)]{Julg-Kasparov-95}}]
\label{theo:Irreducible-decomposition-of-Rumin-complex}
Given $(q,j,i,p)$, we list up all $(s,t)$ such that $s + t \leq n$ and $\Psi_{(q , j , i , p )}^{(s, t)} \not = \{ 0 \}$ as the following:
\begin{enumerate}\itemsep=0pt \setlength{\leftskip}{0.80cm}
\item[Case {\rm I}:]
For $i = j = 0$ and $p = q =0$,
the space is
\[
 \Psi_{(0 , 0 , 0 , 0 )}^{(0, 0)}
 .
\]
\item[Case {\rm II}:]
For $i + j \leq n-2$, $p \geq 1$ and $q \geq 1$,
the spaces are
\[
 \Psi_{(q , j , i , p )}^{(i, j)},
 \quad
 \Psi_{(q , j , i , p )}^{(i+1, j)},
 \quad
 \Psi_{(q , j , i , p )}^{(i, j+1)},
 \quad
 \Psi_{(q , j , i , p )}^{(i+1, j+1)}.
\]
\item[Case {\rm III}:]
For $0 \leq i \leq n-1$, $j=0$, $p \geq 1 $ and $q = 0$,
the spaces are
\[
 \Psi_{(0 , 0 , i , p )}^{(i, 0)},
 \quad
 \Psi_{(0 , 0 , i , p )}^{(i+1, 0)}.
\]
\item[Case {\rm IV}:]
For $i=0$, $0 \leq j \leq n-1$, $p = 0$ and $q \geq 1 $,
the spaces are
\[
 \Psi_{(q , j , 0 , 0 )}^{(0, j)},
 \quad
 \Psi_{(q , j , 0 , 0 )}^{(0, j+1)}.
\]
\item[Case {\rm V}:]
For $i + j = n-1$, $p \geq 1$ and $q \geq 1$,
the spaces are
\[
 \Psi_{(q , j , i , p )}^{(i, j)},
 \quad
 \Psi_{(q , j , i , p )}^{(i+1, j)},
 \quad
 \Psi_{(q , j , i , p )}^{(i, j+1)}.
\]
\item[Case {\rm VI}:]
$i = n-1$, $j=0$, $p \geq 1$ and $q = -1$,
the space is
\[
 \Psi_{(-1 , 0 , n-1 , p )}^{(n, 0)}.
\]
\item[Case {\rm VII}:]
$i=0$, $j=n-1$, $p=-1$ and $q \geq 1$,
the space is
\[
 \Psi_{(q , n-1 , 0 , -1 )}^{(0, n)}.
\]
\end{enumerate}
\end{Proposition}

\begin{Remark}
 About Case VI,
 we substitute $s=n$, $t=0$, $q=-1$, $j=0$, $i=n-1$, for $\Psi_{(q , j , i , p )}^{(s, t)}$,
\[
 \Psi_{(-1 , 0 , n-1 , p )}^{(n, 0)}
 = \cE ^{n,0}
 \cap V (-1, \underline{1}_{0}, \underline{0}_{0}, \underline{-1}_{n-1}, -p )
 = \cE ^{n,0}
 \cap V ( \underline{-1}_{n}, -p )
 .
\]
As the same way, about VII we obtain
\[
 \Psi_{(q , n-1 , 0 , -1 )}^{(0, n)}
 = \cE ^{0,n}
 \cap V (q ,\underline{1}_{n-1}, \underline{0}_{0}, \underline{-1}_{0}, -(-1) )
 = \cE ^{0,n}
 \cap V ( q, \underline{1}_{n} )
 .
\]

\end{Remark}

In~\cite{Kitaoka-19},
$\Delta_{\cE}$ acts as a scalar on these spaces,
and
eigenvalues of $\Delta_\RN $
on the standard CR spheres are explicitly determined.
\begin{Proposition}[{\cite[Theorem 0.1]{Kitaoka-19}}]
\label{prop:all eigenvalues of the Rumin Laplacian on sphere}
Let $S$ be the standard CR sphere with the contact form $\theta = \sqrt{-1} \big(\bar{\partial} - \partial\big) |z|^2$.
Then, on the subspaces of the complexification of $\cE ^\bullet (S)$
corresponding to the representations
$\Psi_{(q,j,i,p)}^{(\bullet , \bullet )}$,
the eigenvalue of $\Delta_\RN $ is
\[
\frac{\big( ( p +i )(q +n -i) + ( q +j )(p +n -j) \big)^2}{4 (n-i-j)^2 } .
\]
\end{Proposition}

\section{Contact analytic torsion of flat vector bundles}
\label{sec:on_flat_vector_bundles}

Let $\mu , \nu_1 , \ldots , \nu_{n+1}$ be integers such that the $\nu_j$ are coprime to $\mu$.
Let $\Gamma$ be the subgroup of~$\big(S^1\big)^{n+1}$ generated by
\begin{equation}
 \gamma
 = ( \gamma_1 , \ldots , \gamma_{n+1} )
 := \big( \exp \big( 2 \pi \ii \nu_1 / \mu \big), \dots , \exp \big( 2 \pi \ii \nu_{n+1} / \mu \big) \big) .
 \notag
\end{equation}
We denote the lens space by
\[
 K
 :=
 S^{2n+1} / \Gamma
 .
\]

Fix $u \in \{ 1, \ldots , \mu \}$ and consider the unitary representation $\alpha_u \colon \pi_1 (K) = \Gamma \to \mathrm{U} (1)$, defined by
\[
 \alpha_u \big( \gamma^{\ell} \big) := \exp \big( 2 \pi \ii u \ell / \mu \big)
 \qquad \text{for} \ \ell \in \mathbb{Z}
 ,
\]
where
\begin{equation}
 \gamma^{\ell}
 := \big( \gamma_1^{\ell} , \ldots , \gamma_{n+1}^{\ell} \big)
 = \big( \exp \big( 2 \pi \ii \nu_1 \ell / \mu \big), \dots , \exp \big( 2 \pi \ii \nu_{n+1} \ell / \mu \big) \big) .
 \notag
\end{equation}
Let $E_{\alpha}$ be the flat vector bundle associated with the unitary representation ${\alpha \colon \pi_1 (K) = \Gamma \! \to \! \mathrm{U} (r)}$, and $E_{u} := E_{\alpha_u}$, which can be considered as $\alpha_u$-equivariant functions on $S^{2n+1}$.

For each unitary representation $(V , \rho )$ of $\mathrm{U}(n+1)$, we define the vector subspace $V^{\alpha_u}$ of $V$ by
\[
V^{\alpha_u} := \big\{ \phi \in V \, \big| \, \rho ( \gamma ) \phi = \alpha ( \gamma )^{-1} \phi \big\}
.
\]

\begin{Proposition}
 We have
 \begin{equation} \label{eq:decom_of_kappa}
 \kappa_\RN (K, E_{u} , g_{\theta , J} ) (s)
 =
 \kappa_{1} (K, E_{u} , g_{\theta , J} ) (s)
 + \kappa_{2} (K, E_{u} , g_{\theta , J} ) (s)
 + \kappa_{3} (K, E_{u} , g_{\theta , J} ) (s),
 \end{equation}
 where
 \begin{align}
 & \kappa_{1} (K, E_{u} , g_{\theta , J} ) (s)
 := -(n+1 ) \dim V^{\alpha_u} ( \underline{0}_{n+1} )
 = \begin{cases}
 -(n+1 ), & u = 0,
 \\
 0, & u \not = 0,
 \end{cases}
 \label{eq:def_of_kappa_1}
 \\
 & \kappa_{2} (K, E_{u} , g_{\theta , J} ) (s)
 :=
 (-1 )^{1} (n+1)
 \sum_{q \geq 1 }
 \frac{
 \dim V^{\alpha_u} ( q, \underline{0}_{n} )
 }{
 \bigl( q /2 \bigr)^{2s}
 }
 \notag
 \\
 & \qquad{}+
 \sum_{j=1}^{n}
 (-1 )^{j+1} (n+1-j)
 \sum_{q \geq 1 }
 \left(
 \frac{
 \dim V^{\alpha_u} ( q, \underline{1}_{j} ,\underline{0}_{n-j} )
 }{
 \bigl( ( q+j)/2 \bigr)^{2s}
 }
 +
 \frac{
 \dim V^{\alpha_u} ( q, \underline{1}_{j-1} ,\underline{0}_{n-j+1} )
 }{
 \bigl( (q+j-1)/2 \bigr)^{2s}
 }
 \right) ,\nonumber
\\
& \kappa_{3} (K, E_{u} , g_{\theta , J} ) (s)
 :=
 \kappa_{2} (K, E_{-u} , g_{\theta , J} ) (s)
 ,
 \label{eq:def_of_kappa_3}
 \end{align}
\end{Proposition}

\begin{proof}
 From Proposition
 \ref{prop:all eigenvalues of the Rumin Laplacian on sphere}, we check that
 the terms of $\kappa_{\RN} ( K, E_{u} , g ) (s)$ in Cases II and V in Proposition~\ref{theo:Irreducible-decomposition-of-Rumin-complex} cancel each other.
 ``The sum of the terms of $\kappa_\RN (K, E_{u} , g_{\theta , J}) (s)$ in Case II'' is
 \begin{align*}
 &
 \sum_{a=0}^{n-2}
 \sum_{i+j=a}
 \sum_{\substack{p\geq 1,\\ q\geq 1}}
 \big( (-1)^{i+j+1} (n+1 -i-j)
 \\
 & \qquad\quad{}
 + 2 (-1)^{i+j+2} (n -i-j) + (-1)^{i+j+3} (n-1 -i-j) \big)
 \\
 & \qquad\quad{}\times
 \frac{
 \dim V^{\alpha_u} ( q, \underline{1}_{j}, \underline{0}_{n-1-i-j}, \underline{-1}_{i}, -p )
 }{
 \big( ((p+i)(q+n-i) + (q+j)(p+n-j)) / 2 (n-i-j) \big)^{2s}
 }
 \\
 &\qquad{}=
 \sum_{a=0}^{n-2}
 \sum_{i+j=a}
 \sum_{\substack{p\geq 1,\\ q\geq 1}}
 \bigl( 0 \bigr)
 \frac{
 \dim V^{\alpha_u} ( q, \underline{1}_{j}, \underline{0}_{n-1-i-j}, \underline{-1}_{i}, -p )
 }{
 \big( ((p+i)(q+n-i) + (q+j)(p+n-j)) / 2 (n-i-j) \big)^{2s}
 } = 0.
 \end{align*}
 Similarly,
 \text{``the sum of the terms of $\kappa_\RN (K, E_{u} , g_{\theta , J}) (s)$ in Case V''} is
 \begin{gather*}
 \sum_{i+j=n-1}
 \sum_{\substack{p\geq 1,\\ q\geq 1}}
 \big( (-1)^{n} (n+1 -(n-1)))
 + 2 (-1)^{n+1} (n+1 -n) \big)
 \\
  \qquad\quad{}\times
 \frac{
 \dim V^{\alpha_u} ( q, \underline{1}_{j}, \underline{-1}_{i}, -p )
 }{
 \big( ((p+i)(q+n-i) + (q+j)(p+n-j)) / 2 (n+1- i-j) \big)^{2s}
 }
 \\
\qquad{}=
 \sum_{i+j=n-1}
 \sum_{\substack{p\geq 1, \\ q\geq 1}}
 \bigl( 0 \bigr)
 \frac{
 \dim V^{\alpha_u} ( q, \underline{1}_{j}, \underline{-1}_{i}, -p )
 }{
 \big( ((p+i)(q+n-i) + (q+j)(p+n-j)) / 2 (n-i-j) \big)^{2s}
 }
 = 0.
 \end{gather*}

 The function $\kappa_{1} (K, E_{u} , g_{\theta , J} ) (s)$
 is the sum of the terms of $\kappa_\RN (K, E_{u} , g_{\theta , J}) (s)$ in Case I.

 Next we consider the sum of the terms of $\kappa_{\cE}$ in Cases III and VI.
 For $j$,
 ``the sum of the terms of $\kappa_{\cE} (K, E_{u} , g_{\theta , J} ) $ in $\cE^{0,j}$, Cases III and VI'' is
 \begin{gather*}
 (-1 )^{1} (n+1)
 \sum_{q \geq 1 }
 \frac{
 \dim V^{\alpha_u} ( q ,\underline{0}_{n} )
 }{
 \bigl( q/2 \bigr)^{2s}
 }
 ,
\quad
 j = 0,
 \\
 (-1 )^{j+1} (n+1-j)
 \sum_{q \geq 1 }
 \left(
 \frac{
 \dim V^{\alpha_u} ( q, \underline{1}_{j} ,\underline{0}_{n-j} )
 }{
 \bigl( ( q+j)/2 \bigr)^{2s}
 }
 +
 \frac{
 \dim V^{\alpha_u} ( q, \underline{1}_{j-1} ,\underline{0}_{n-j+1} )
 }{
 \bigl( (q+j-1)/2 \bigr)^{2s}
 }
 \right)
 ,\quad
 1 \leq j \leq n.
 \end{gather*}
 ``The sum of the terms of $\kappa_{\cE}$ in Cases III and VI'' is
 \begin{align}
 &
 (-1 )^{1} (n+1)
 \sum_{q \geq 1 }
 \frac{
 \dim V^{\alpha_u} ( q, \underline{0}_{n} )
 }{
 \bigl( q /2 \bigr)^{2s}
 }
 \notag
 \\
 & \qquad\quad{}+
 \sum_{j=1}^{n}
 (-1 )^{j+1} (n+1-j)
 \sum_{q \geq 1 }
 \left(
 \frac{
 \dim V^{\alpha_u} ( q, \underline{1}_{j} ,\underline{0}_{n-j} )
 }{
 \bigl( ( q+j)/2 \bigr)^{2s}
 }
 +
 \frac{
 \dim V^{\alpha_u} ( q, \underline{1}_{j-1} ,\underline{0}_{n-j+1} )
 }{
 \bigl( (q+j-1)/2 \bigr)^{2s}
 }
 \right)
 \notag
 \\
 &\qquad{}= \kappa_{2} (K, E_{u} , g_{\theta , J} ) (s) .
 \label{eq:kappa2_formula-01}
 \end{align}

 Finally, we consider the sum of the terms of $\kappa_{\cE}$ in Cases IV and VII.
 As the same way \eqref{eq:kappa2_formula-01} in Cases III and VI,
 ``the sum of the terms of $\kappa_{\cE}$ in Cases IV and VII''
 is given by
 \begin{align}
 &
 (-1 )^{1} (n+1)
 \sum_{q \geq 1 }
 \frac{
 \dim V^{\alpha_u} ( \underline{0}_{n}, -q )
 }{
 \bigl( q /2 \bigr)^{2s}
 }
 \label{eq:kappa3_formula-01}
 \\
 &{}+
 \sum_{j=1}^{n}
 (-1 )^{j+1} (n+1-j)
 \sum_{q \geq 1 }
 \left(
 \frac{
 \dim V^{\alpha_u} ( \underline{0}_{n-j}, \underline{-1}_{j} , -q )
 }{
 \bigl( ( q+j)/2 \bigr)^{2s}
 }
 +
 \frac{
 \dim V^{\alpha_u} ( \underline{0}_{n-j+1}, \underline{-1}_{j-1} , -q )
 }{
 \bigl( (q+j-1)/2 \bigr)^{2s}
 }
 \right). \notag
 \end{align}
 Let $(V , \rho)$ be the unitary representation of $\mathrm{U}(n+1)$.
 We define the representation $(\overline{V}, \overline{\rho})$ by
 \[
 \overline{V} := V,
 \qquad
 \overline{\rho} (U)
 :=
 \overline{\rho (U)}
 \qquad \text{for} \ U \in \mathrm{U} (n+1)
 .
 \]
 Since $(V ( q, \underline{1}_{j} ,\underline{0}_{n-j} ) ,\rho )$ is the unitary representation,
 its conjugate representation is isomorphic to its dual representation as $\mathrm{U}(n+1)$-module.
 From~\cite[Theorem 3.2.13]{Goodman-Wallach-09},
 the conjugate representation of $( V ( q, \underline{1}_{j} ,\underline{0}_{n-j} ) , \rho )$ is isomorphic to
 $(V ( \underline{0}_{n-j} , \underline{-1}_{j} , -q ) , \overline{\rho} )$ as $\mathrm{U}(n+1)$-module.
 Therefore, we have
 \begin{align*}
 V^{\alpha_{-u}} ( q, \underline{1}_{j} ,\underline{0}_{n-j} )
 & =
 \big\{ \phi \in V ( q, \underline{1}_{j} ,\underline{0}_{n-j}) \, \big| \, \alpha_{-u} (- \gamma ) \phi = \rho ( \gamma ) \phi \big\}
 \\
 & =
 \big\{ \phi \in V ( q, \underline{1}_{j} ,\underline{0}_{n-j}) \, \big| \, \alpha_u (- \gamma )^{-1} \phi = \rho ( \gamma ) \phi \big\}
 \\
 & =
 \big\{ \phi \in V ( q, \underline{1}_{j} ,\underline{0}_{n-j}) \, \big| \, \alpha_u ( \gamma ) \phi = \rho ( \gamma ) \phi \big\}
 \\
 & \cong
 \big\{ \phi \in V ( \underline{0}_{n-j}, \underline{-1}_{j} , -q ) \, \big| \, \overline{\alpha_u ( \gamma )} \phi = \overline{\rho} ( \gamma ) \phi \big\}
 \\
 & =
 \big\{ \phi \in V ( \underline{0}_{n-j}, \underline{-1}_{j} ,-q ) \, \big| \, \alpha_u ( \gamma )^{-1} \phi = \overline{\rho} ( \gamma ) \phi \big\}
 \\
 & = V ( \underline{0}_{n-j} , \underline{-1}_{j} , -q )^{\alpha_{u}}
 ,
 \end{align*}
 where $\cong$ means isomorphic as real vector spaces via the complex conjugate.
 Then, from \eqref{eq:kappa2_formula-01} and
 \eqref{eq:kappa3_formula-01},
 ``the sum of the terms of $\kappa_{\cE}$ in Cases IV and VII''
 is given by
 \begin{align*}
 &
 (-1 )^{1} (n+1)
 \sum_{q \geq 1 }
 \frac{
 \dim V^{\alpha_u} ( \underline{0}_{n}, -q )
 }{
 \bigl( q /2 \bigr)^{2s}
 }
 \notag
 \\
 & \quad +
 \sum_{j=1}^{n}
 (-1 )^{j+1} (n+1-j)\!
 \sum_{q \geq 1 }\!
 \left(
 \frac{
 \dim V^{\alpha_u} ( \underline{0}_{n-j}, \underline{-1}_{j} , -q )
 }{
 \bigl( ( q+j)/2 \bigr)^{2s}
 }
 +
 \frac{
 \dim V^{\alpha_u} ( \underline{0}_{n-j+1}, \underline{-1}_{j-1} , -q )
 }{
 \bigl( (q+j-1)/2 \bigr)^{2s}
 }
 \right)
 \notag
 \\
 & =
 (-1 )^{1} (n+1)
 \sum_{q \geq 1 }
 \frac{
 \dim V^{\alpha_{-u}} ( q, \underline{0}_{n} )
 }{
 \bigl( q /2 \bigr)^{2s}
 }
 \notag
 \\
 & \quad +
 \sum_{j=1}^{n}
 (-1 )^{j+1} (n+1-j)
 \sum_{q \geq 1 }
 \left(
 \frac{
 \dim V^{\alpha_{-u}} ( q, \underline{1}_{j} ,\underline{0}_{n-j} )
 }{
 \bigl( ( q+j)/2 \bigr)^{2s}
 }
 +
 \frac{
 \dim V^{\alpha_{-u}} ( q, \underline{1}_{j-1} ,\underline{0}_{n-j+1} )
 }{
 \bigl( (q+j-1)/2 \bigr)^{2s}
 }
 \right)
 \notag
 \\
 & =
 \kappa_{2} (K, E_{-u} , g_{\theta , J} ) (s).
 \tag*{\qed}
 \end{align*}
\renewcommand{\qed}{}
\end{proof}

We set for $q \geq 1$,
\[
 V ( q, \underline{1}_{-1} ,\underline{0}_{n+1} )
 := \{ 0 \}
 .
\]
We have
\begin{align}
 & \kappa_{2} (K, E_{u} , g_{\theta , J} ) (s)
 =
 (-1 )^{1} (n+1)
 \sum_{q \geq 1 }
 \frac{
 \dim V^{\alpha_u} ( q, \underline{0}_{n} )
 }{
 \bigl( q /2 \bigr)^{2s}
 }
 \notag
 \\
 & \qquad\quad{}+
 \sum_{j=1}^{n}
 (-1 )^{j+1} (n+1-j)
 \sum_{q \geq 1 }
 \left(
 \frac{
 \dim V^{\alpha_u} ( q, \underline{1}_{j} ,\underline{0}_{n-j} )
 }{
 \bigl( ( q+j)/2 \bigr)^{2s}
 }
 +
 \frac{
 \dim V^{\alpha_u} ( q, \underline{1}_{j-1} ,\underline{0}_{n-j+1} )
 }{
 \bigl( (q+j-1)/2 \bigr)^{2s}
 }
 \right)
 \notag
 \\
 &\qquad{}=
 \sum_{j=0}^{n}
 (-1 )^{j+1} (n+1-j)
 \sum_{q \geq 1 }
 \left(
 \frac{
 \dim V^{\alpha_u} ( q, \underline{1}_{j} ,\underline{0}_{n-j} )
 }{
 \bigl( ( q+j)/2 \bigr)^{2s}
 }
 +
 \frac{
 \dim V^{\alpha_u} ( q, \underline{1}_{j-1} ,\underline{0}_{n-j+1} )
 }{
 \bigl( (q+j-1)/2 \bigr)^{2s}
 }
 \right)
 \notag
 \\
 &\qquad{}=
 \sum_{j=0}^{n}
 (-1 )^{j+1} (n+1-j)
 \sum_{q \geq 1 }
 \frac{
 \dim V^{\alpha_u} ( q, \underline{1}_{j} ,\underline{0}_{n-j} )
 +
 \dim V^{\alpha_u} ( q+1, \underline{1}_{j-1} ,\underline{0}_{n-j+1} )
 }{
 \bigl( (q+j)/2 \bigr)^{2s}
 }
 \notag
 \\
 & \qquad\quad{}
 +
 \sum_{j=1}^{n}
 (-1 )^{j+1} (n+1-j)
 \frac{
 \dim V^{\alpha_u} ( \underline{1}_{j} ,\underline{0}_{n-j+1} )
 }{
 \bigl( j/2 \bigr)^{2s}
 }
 .
 \notag
\end{align}

Let $\chi_V$ be the character of the representation $(V, \rho)$ of $ \mathrm{U} (n+1)$.
We note that for each representation $(V,\rho)$,
\[
 \dim V^{\alpha_u}
 = \sum_{t \in \Gamma } \chi_{V} ( t ) \alpha_u ( t ) / \# \Gamma
 ,
\]
(cf.~\cite[equation~(2.9)]{Fulton-Harris-04}).
By Littlewood--Richardson's rule (cf.~\cite[Corollary~3]{Fulton-97}), we have
\begin{equation}
 \chi_{V ( \underline{1}_{j} ,\underline{0}_{n-j+1} )}
 \chi_{V ( q, \underline{0}_{n} )}
 =
 \chi_{V ( q, \underline{1}_{j} ,\underline{0}_{n-j} )}
 +
 \chi_{V ( q+1, \underline{1}_{j-1} ,\underline{0}_{n-j+1} )}
 .
 \label{eq:Rich-Litt}
\end{equation}

From \eqref{eq:Rich-Litt},
\begin{align}
 & \kappa_{2} (K, E_{u} , g_{\theta , J} ) (s)
 =
 \frac{1}{\mu}
 \sum_{j=0}^{n}
 (-1 )^{j+1} (n+1-j)
 \sum_{q \geq 1 }
 \sum_{\ell = 0}^{\mu -1}
 \frac{
 \chi_{V ( \underline{1}_{j} ,\underline{0}_{n-j+1} )} \big( \gamma^{\ell} \big)
 \chi_{V ( q, \underline{0}_{n} )} \big( \gamma^{\ell} \big)
 \alpha_u \big( \gamma^{\ell} \big)
 }{
 \bigl( (q+j)/2 \bigr)^{2s}
 }
 \notag
 \\
 & \qquad\quad{}
 +
 \frac{1}{\mu}
 \sum_{j=1}^{n}
 (-1 )^{j+1} (n+1-j)
 \sum_{\ell = 0}^{\mu -1}
 \frac{
 \chi_{V ( \underline{1}_{j} ,\underline{0}_{n-j+1} )} \big( \gamma^{\ell} \big)
 \alpha_u \big( \gamma^{\ell} \big)
 }{
 \bigl( j/2 \bigr)^{2s}
 }
 \notag
 \\
 & \qquad{}=
 \frac{
 2^{2s}
 }{
 \mu \Gamma (2s)
 }
 \sum_{j=0}^{n}
 (-1 )^{j+1} (n+1-j)
 \notag
 \\
 & \qquad\quad\quad{}\times
 \sum_{q \geq 1 }
 \sum_{\ell = 0}^{\mu -1}
 \int_{0}^{\infty}
 \chi_{V ( \underline{1}_{j} ,\underline{0}_{n-j+1} )} \big( \gamma^{\ell} \big)
 \chi_{V ( q, \underline{0}_{n} )} \big( \gamma^{\ell} \big)
 \alpha_u \big( \gamma^{\ell} \big)
 {\rm e}^{-(j + q) x }
 x^{2s-1}
 \, {\rm d}x
 \notag
 \\
 & \qquad\quad{}
 +
 \frac{
 2^{2s}
 }{
 \mu \Gamma (2s)
 }
 \sum_{j=1}^{n}
 (-1 )^{j+1} (n+1-j)
 \sum_{\ell = 0}^{\mu-1}
 \int_0^{\infty}
 \chi_{V ( \underline{1}_{j} ,\underline{0}_{n-j+1} )} \big( \gamma^{\ell} \big)
 \alpha_u \big( \gamma^{\ell} \big)
 {\rm e}^{-jx}
 x^{2s-1} \, {\rm d}x
 \notag
 \\
 & \qquad{}=
 \frac{
 2^{2s}
 }{
 \mu \Gamma (2s)
 }
 \sum_{\ell = 0}^{\mu-1}
 \int_{0}^{\infty}
 \Biggl(
 \sum_{j=0}^{n}
 (-1 )^{j+1} (n+1-j)
 \chi_{V ( \underline{1}_{j} ,\underline{0}_{n-j+1} )} \big( \gamma^{\ell} \big)
 {\rm e}^{-j x }
 \sum_{q \geq 1 }
 \chi_{V ( q, \underline{0}_{n} )} \big( \gamma^{\ell} \big)
 {\rm e}^{-q x }
 \notag \\
 & \qquad\quad{}
 +
 \sum_{j=1}^{n}
 (-1 )^{j+1} (n+1-j)
 \chi_{V ( \underline{1}_{j} ,\underline{0}_{n-j+1} )} \big( \gamma^{\ell} \big)
 {\rm e}^{-jx}
 \Biggr)
 \alpha_u \big( \gamma^{\ell} \big)
 x^{2s-1} \, {\rm d}x
 \label{eq:kappa_heatkernel_rep}
 .
\end{align}
We consider the contents of the integral for the last equation of \eqref{eq:kappa_heatkernel_rep}.
It is known that for $t = (t_1 , \dots , t_{n+1} ) \in \big(S^1\big)^{n+1}$,
\begin{align}
& \chi_{V ( \underline{1}_{j} ,\underline{0}_{n-j+1} )} ( t )
 =
 \sum_{\substack{\beta_1 + \cdots + \beta_{n+1} = j , \\
 0 \leq \beta_1 , \ldots , \beta_{n+1} \leq 1}}
 t_{1}^{\beta_1} \cdots t_{n+1}^{\beta_{n+1}}
 \label{eq:chi_j}
 ,
 \\
& \chi_{V ( q, \underline{0}_{n} )} ( t )
 =
 \sum_{\substack{ \alpha_1 + \cdots + \alpha_{n+1} = q\\
 \alpha_1 , \ldots , \alpha_{n+1} \geq 0}}
 t_{1}^{\alpha_1} \cdots t_{n+1}^{\alpha_{n+1}}
 \label{eq:chi_q}
 ,
\end{align}
(cf.~\cite[equations (6.1) and (6.2)]{Fulton-Harris-04}).
We set $X := {\rm e}^{-x}$
and
\[
 F_1 ( t , X)
 :=
 \sum_{j=0}^{n+1}
 (-1)^j
 \chi_{V ( \underline{1}_{j} ,\underline{0}_{n-j+1} )} ( t )
 X^j
 .
\]
Then \eqref{eq:chi_j} gives
\begin{align}
& F_1 ( t ,X)
 =
 \prod_{j=1}^{n+1}
 (1 - t_j X)
 ,
 \label{eq:F_1_formula}
 \\
& X \frac{\partial F_1}{\partial X} ( t , X)
 =
 \sum_{j=0}^{n+1}
 (-1)^j
 j
 \chi_{V ( \underline{1}_{j} ,\underline{0}_{n-j+1} )} ( t )
 X^j
 =
 -
 \sum_{i=1}^{n+1}
 t_i X
 \prod_{j=1, j \not = i }^{n+1}
 (1 - t_j X)
 .
 \label{eq:F_1_derivative}
\end{align}
From the definition of $F_1$ and \eqref{eq:F_1_derivative},
we have
\begin{align}
 \sum_{j=0}^{n}
 (-1 )^{j+1} (n+1-j)
 \chi_{V ( \underline{1}_{j} ,\underline{0}_{n-j+1} )} ( t )
 X^j
 =
 - (n+1 ) F_1 ( t , X )
 + X \frac{\partial F_1 }{\partial X} ( t , X )
 .
 \label{eq:kappa_j_term}
\end{align}
We set
\[
 F_2 ( t ,X )
 :=
 \sum_{q \geq 1 }
 \chi_{V ( q, \underline{0}_{n} )} ( t ) X^q .
\]
From \eqref{eq:chi_q} and \eqref{eq:F_1_formula},
we can rewrite $F_2$ as
\begin{equation}
 F_2 ( t , X)
 =
 \prod_{j=1}^{n+1}
 \frac{1}{1 - t_j X} -1
 =
 \frac{1}{F_1 ( t , X)} -1
 .
 \label{eq:kappa_q_term}
\end{equation}
From \eqref{eq:F_1_formula}--\eqref{eq:kappa_q_term},
we can deduce that
\begin{align}
 &
 \sum_{j=0}^{n}
 (-1 )^{j+1} (n+1-j)
 \chi_{V ( \underline{1}_{j} ,\underline{0}_{n-j+1} )} ( t )
 X^j
 \sum_{q \geq 1 }
 \chi_{V ( q, \underline{0}_{n} )} ( t ) X^q
 \notag
 \\
 &
 \qquad \quad{}+
 \sum_{j=1}^{n}
 (-1 )^{j+1} (n+1-j)
 \chi_{V ( \underline{1}_{j} ,\underline{0}_{n-j+1} )} ( t )
 X^j
 \notag
 \\
 & \qquad{}=
 -
 \left(
 (n+1) F_1 ( t,X) - X \frac{\partial F_1}{\partial X} ( t , X )
 \right)
 \left(
 \frac{1}{F_1 ( t , X ) } -1
 \right)
 \notag
 \\
 & \qquad \quad{}
 -
 \left(
 (n+1)( F_1 ( t , X ) -1 )
 - X \frac{\partial F_1}{\partial X} ( t , X )
 \right)
 \notag \\
 & \qquad{}=
 \frac{X \frac{\partial F_1}{\partial X} ( t , X )}{F_1 ( t , X )}
 =
 -
 \sum_{j=1}^{n+1} \frac{t_j X}{1 - t_j X }
 .
 \label{eq:kappa_contents_of_inte}
\end{align}

From \eqref{eq:kappa_heatkernel_rep} and \eqref{eq:kappa_contents_of_inte},
we see
\begin{align}
 \kappa_{2} (K, E_{u} , g_{\theta , J} ) (s)
 & =
 -
 \frac{
 2^{2s}
 }{
 \mu \Gamma (2s)
 }
 \sum_{\ell = 0}^{\mu-1}
 \sum_{j=1}^{n+1}
 \int_{0}^{\infty} \frac{\gamma_j^{\ell} {\rm e}^{-x}}{1 - \gamma_j^{\ell} {\rm e}^{-x}}
 {\rm e}^{ 2 \pi \ii u \ell / \mu } x^{2s-1} \, {\rm d}x
 \notag
 \\
 %
 & =
 -
 \sum_{j=1}^{n+1}
 \frac{
 2^{2s}
 }{
 \mu \Gamma (2s)
 }
 \sum_{\ell = 0}^{\mu -1}
 \int_{0}^{\infty} \sum_{q=1}^{\infty} \exp \big( 2 \pi \ii (q \nu_j + u ) \ell / \mu \big) {\rm e}^{-qx} x^{2s-1} \, {\rm d}x
 \notag \\
 & =
 -
 2^{2s}
 \sum_{j=1}^{n+1}
 \sum_{\ell = 0}^{\mu -1}
 \sum_{q=1}^{\infty}
 \frac{
 \exp \big( 2 \pi \ii (q \nu_j + u ) \ell / \mu \big)
 }{
 \mu
 }
 q^{-2s}
 .
 \label{eq:kappa_2_exp_sum}
\end{align}
Let ${\tau}_j$ be the integers in $\{ 1 , \ldots , \mu \}$ such that ${\tau}_j \nu_j \equiv 1 \mod \mu$.
Since the multiplication of $\nu_j \in ( \mathbb{Z} / \mu \mathbb{Z} )^{\times}$ induced the bijective map from $\mathbb{Z} / \mu \mathbb{Z} $ to $\mathbb{Z} / \mu \mathbb{Z} $,
We have
\begin{align*}
 \sum_{\ell = 0}^{\mu -1}
 \exp \big( 2 \pi \ii (q \nu_j + u ) \ell / \mu \big)
 & =
 \sum_{\ell = 0}^{\mu -1}
 \exp \big( 2 \pi \ii (q + u \tau_j ) \nu_j \ell / \mu \big)
 \\
 & =
 \sum_{\ell = 0}^{\mu -1}
 \exp \big( 2 \pi \ii (q + u \tau_j ) \ell / \mu \big)
 \\
 & =
 \begin{cases}
 0 , & q \not \equiv - u \tau_j \mod \mu ,
 \\
 \mu , & q \equiv - u \tau_j \mod \mu .
 \end{cases}
\end{align*}
For $w \in \mathbb{Z}$ let $A_{\mu} ( w )$ be the integer between $1$ and $\mu$ which is congruent to $w$ modulo $\mu$,
then from \eqref{eq:kappa_2_exp_sum}, we can rewrite $\kappa_{2}$ as
\begin{align}
 \kappa_{2} (K, E_{u} , g_{\theta , J} ) (s)
 & =
 -
 2^{2s}
 \sum_{j=1}^{n+1}
 \sum_{\substack{q >0 , \\ q \equiv - u {\tau}_j \mod \mu } }^{\infty}
 q^{-2s}
 =
 -
 2^{2s}
 \sum_{j=1}^{n+1}
 \sum_{q=0}^{\infty}
 \bigl(
 q \mu + A_{\mu} (- u {\tau}_j)
 \bigr)^{-2s}
 \notag
 \\
 & =
 -
 2^{2s} \mu^{-2s}
 \sum_{j=1}^{n+1}
 \zeta \bigl( 2s , A_{\mu} (-u {\tau}_j) / \mu \bigr)
 ,
 \label{eq:kappa_2}
\end{align}
where for $0 < a \leq 1$, $\zeta (s , a) := \sum_{q =0}^{\infty} (q + a)^{-s}$ is the Hurwitz zeta function.

Next, we calculate $\kappa_{3}$.
As the same way in calculating $\kappa_{2}$,
from \eqref{eq:def_of_kappa_3},
we can rewrite $\kappa_{3}$ as
\begin{equation}
 \kappa_{3} ( K, E_{u} , g_{\theta , J} ) (s)
 =
 -
 2^{2s} \mu^{-2s}
 \sum_{j=1}^{n+1}
 \zeta \bigl( 2s , A_{\mu} ( u {\tau}_j) / \mu \bigr)
 .
 \label{eq:kappa_3}
\end{equation}

From \eqref{eq:decom_of_kappa}, \eqref{eq:def_of_kappa_1}, \eqref{eq:kappa_2} and \eqref{eq:kappa_3},
we have
\[
 \kappa_{\RN} ( K, E_{u} , g_{\theta , J} ) (s)
 =
 \begin{cases}
 \displaystyle
 -(n+1)
 \bigl( 1+
 2^{2s+1}
 \mu^{-2s}
 \zeta (2s)
 \bigr)
 ,
 \quad
 u = 0 .
 \\
 \displaystyle
 -
 2^{2s} \mu^{-2s}
 \sum_{j=1}^{n+1}
 \big(
 \zeta \bigl( 2s , A_{\mu} (u {\tau}_j) / \mu \bigr)
 +
 \zeta \bigl( 2s , 1 - A_{\mu} (u {\tau}_j) / \mu \bigr)
 \big)
 ,
 \quad
 u \not = 0 .
 \end{cases}
\]
It is known that
$\zeta (0) = -1/2$ and $\zeta ' (0) = - \log (2 \pi)/2$
and for $0 < a < 1$,
\begin{align*}
& \zeta (0 ,a) + \zeta (0 , 1-a)
 = 0 , \\
& \zeta ' (0 ,a) + \zeta ' (0 , 1-a)
 = - \log \big| {\rm e}^{2 \pi \ii a} -1 \big|
 .
\end{align*}
Using the above equations,
we conclude
for $u = 0$
\begin{align*}
& \kappa_{\RN} (K, \underline{\bC} , g_{\theta , J} ) (0)
 =
 - (n+1) \bigl( 1 + 2 \zeta (0) \bigl) = 0
 ,
 \\
& \kappa_{\RN} (K, \underline{\bC} , g_{\theta , J} ) ' (0)
 =
 2 (n+1) \log \frac{4 \pi}{\mu}
 ,
\end{align*}
and
for $u \not = 0$
\begin{align*}
& \kappa_\RN (K, E_{u} , g_{\theta , J} ) (0)
 = 0
 ,
 \\
& \kappa_\RN (K, E_{u} , g_{\theta , J} ) ' (0)
 =
 2 \sum_{j=1}^{n+1} \log \big| {\rm e}^{ 2 \pi \ii u {\tau}_j / \mu } -1 \big|
\end{align*}
as claimed.

\section{Ray--Singer torsion of the trivial bundle}
\label{sec:Ray--Singer_torsion_on_the_trivial_vector_bundle}

We compute the analytic torsion of the trivial bundle on lens spaces.
We define the Ray--Singer torsion function associated with $\big(\Omega^{\bullet}(M,E) , {\rm d}^{\nabla}\big)$ by
\[
 \kappa_{\deRham} (M, E , g ) (s)
 := \sum_{k=0}^{2n+1} (-1)^k k \zeta \big(\Delta_{\deRham , g}^k \big) (s )
 ,
\]
where $\zeta \big(\Delta_{\deRham , g}^k \big) (s )$ is the spectral zeta function of the $k$-th Hodge--de~Rham Laplacian $\Delta_{\deRham , g}^k$
with respect to $g$.
We define the Ray--Singer torsion $T_{\deRham}$
by
\[
 2 \log T_{\deRham} (M, E , g )
 := \kappa_{\deRham} (M, E , g )' (0)
 .
\]
Following the derivation of~\cite[equation (3)]{Ray-70}, we have
\[
 2 \log T_{\deRham} ( K , \underline{\bC} , g_{\mathrm{std}} )
 =
 \frac{1}{\mu}
 \sum_{\ell = 0}^{\mu -1}
 \sum_{j=0}^{2n}
 (-1)^{j+1}
 \zeta ' \big(0 ; j, \gamma^{\ell} \big)
 ,
\]
where $\lambda_m$ is the $m$-th eigenvalue of ${\rm d}^{\adjoint} {\rm d}$
and
\begin{align*}
& \zeta \big(s ; j, \gamma^{\ell} \big)
 =
 \sum_{m=0}^{\infty}
 \lambda_m^{-s} \Tr \big( \gamma^{\ell} \rest{X_{j,m}}\big)
 ,
 \\
& X_{j,m}
 :=
 \big\{
 \phi \in \Omega^j (K)
 \,|\,
{\rm d}^{\adjoint} {\rm d} \phi = \lambda_m \phi
 \big\}
 .
\end{align*}
We recall from~\cite[p.~123]{Ray-70}
\begin{align*}
 &
 \sum_{j=0}^{2n}
 (-1)^{j+1}
 \zeta \big(s ; j, \gamma^{\ell} \big)
 \\
 & \qquad=
 \big(\Gamma (s)\big)^{-2}
 \int_0^{\infty} t^{2s-1}
 \int_0^1 \bigl( u(1-u) \bigr)^{s-1}
 \big(
 f_{0} (t, u)
 + f_{1} (t , \ell \nu / \mu )
 \big) \, {\rm d}u \, {\rm d}t
 + \mathcal{O} \big(s^2\big)
 ,
\end{align*}
where $\nu = (\nu_1 , \ldots , \nu_{n+1}) $ and for $\sigma = (\sigma_1 , \ldots , \sigma_{n+1} )\in \bR^{n+1}$,
\begin{align*}
 f_{0} (t, u)
 & =
 \sum_{j=0}^{n+1}
 (-1)^{j} \binom{n}{j}
 \left(
 \frac{
 2 \sinh t u \sinh t (1-u)
 }{
 \sinh t
 }
 \right)^{j}
 - \frac{{\rm e}^{-(2n+1)tu}}{2 \sinh tu}
 - \frac{{\rm e}^{-(2n+1)t(1-u)}}{2 \sinh t(1-u)}
 ,
 \\
 f_{1} (t , \sigma )
 & =
 \sum_{k}^{n+1}
 \left(
 1
 -
 \frac{
 \sinh t
 }{
 \cosh t - \cos 2 \pi \sigma_k
 }
 \right)
 .
\end{align*}
We set
\begin{align*}
& h_0 (s) :=
 \big(\Gamma (s)\big)^{-2}
 \int_0^{\infty} t^{2s-1}
 \int_0^1 \bigl( u(1-u) \bigr)^{s-1}
 f_0 (t,u ) \, {\rm d}u \, {\rm d}t
 ,
 \\
& h_1 \big(s , \gamma^{\ell}\big) :=
 \big(\Gamma (s)\big)^{-2}
 \int_0^{\infty} t^{2s-1}
 \int_0^1 \bigl( u(1-u) \bigr)^{s-1}
 f_1 (t, \ell \nu/\mu ) \, {\rm d}u \, {\rm d}t
 .
\end{align*}
From~\cite[p.~125]{Ray-70},
it is seen that
\begin{equation*}
 h_1 \big(s , \gamma^{\ell}\big)
 =
 - \frac{1}{2}
 \sum_{j=0}^{\mu -1}
 \Tr \big(\gamma^{j \ell}\big)
 \bigl(
 \zeta (2s , j /\mu)
 +
 \zeta (2s , 1- j /\mu)
 \bigr)
 -
 2 (n+1) \mu^{-2s} \zeta (2s)
 ,
\end{equation*}
where
\[
 \Tr \big(\gamma^{\ell}\big) =
 \sum_{j=1}^{n+1}
 \big(
 {\rm e}^{2 \pi \ii \ell \nu_j / \mu} + {\rm e}^{-2 \pi \ii \ell \nu_j / \mu}
 \big)
 .
\]
Taking the average of $h_1 (s , \gamma^{\ell})$, we have
\[
 \frac{1}{\mu} \sum_{\ell = 0}^{\mu-1} h_1 \big(s , \gamma^{\ell}\big)
 =
 - 2 (n+1)
 \mu^{-2s} \zeta (2s)
\]
Using $\zeta '(0) = - \log (2 \pi) /2$, we get
\begin{equation}
 \frac{1}{\mu} \sum_{\ell = 0}^{\mu-1} h_1 ' \big(0 , \gamma^{\ell}\big)
 = 2 (n+1) \log \left( \frac{2 \pi}{ \mu} \right)
 .
 \label{eq:f_1_at_origin}
\end{equation}
We recall the Ray--Singer torsion on spheres,
\begin{Proposition}[{\cite{Weng-You-96}}]
 \label{prop:Weng-You}
 \[
 T_{\deRham} (S, \underline{\bC} , g_{\mathrm{std}})
 = \frac{2 \pi^{n+1}}{n!}
 \]
\end{Proposition}
\begin{Remark}
 The Ray--Singer torsion of spheres can easily be determined using the Cheeger--M\"uller theorem, a result which predates~\cite{Weng-You-96}.
 Details can be found in~\cite{Melo-Spreafico-09}.
\end{Remark}
To put $\mu =1$, $\nu = (1 , \ldots , 1)$, from \eqref{eq:f_1_at_origin}, it follows that
\begin{align}
 h_0' (0)
 & =
 2 \log T_{\deRham} (S, \underline{\bC} , g_{\mathrm{std}})
 - \frac{1}{\mu} \sum_{\ell = 0}^{\mu-1} h_1 ' \big(0 , \gamma^{\ell}\big)
 \notag \\
 & =
 2 \log \left( \frac{2 \pi^{n+1}}{ n!} \right)
 - 2 (n+1) \log ( 2 \pi )
 = 2 \log \left( \frac{2^{-n}}{ n!} \right)
 .
 \label{eq:f_0_at_origin}
\end{align}
By \eqref{eq:f_1_at_origin} and \eqref{eq:f_0_at_origin},
we conclude
\begin{align}
 2 \log T_{\deRham} (K, \underline{\bC} , g_{\mathrm{std}})
 & =
 \frac{1}{\mu} \sum_{\ell = 0}^{\mu -1} h_0 ' (0 )
 + \frac{1}{\mu} \sum_{\ell = 0}^{\mu -1} h_1 ' \big(0 , \gamma^{\ell}\big)
 \notag \\
 & = 2 \log \left( \frac{2^{-n}}{n!} \right) + 2 (n+1) \log \left( \frac{2 \pi}{\mu} \right)
 = 2 \log \left( \frac{2 \pi^{n+1}}{ n! \mu^{n+1}} \right)
 .
 \label{eq:pre_Proposition_trivial_bundle}
\end{align}

For the metric $g$ on $E$ over $M$,
we set $g_{\rho} := {\rm e}^{2 \rho} g$ for $\rho \in \mathbb{R}$.
Then, the Hodge--de~Rham Laplacian $\Delta_{\deRham , g_{\rho}}$ with respect to $g_{\rho}$
is given by
\begin{equation}
 \Delta_{\deRham , g_{\rho}} = {\rm e}^{ -2 \rho } \Delta_{\deRham , g}
 \label{eq:scaling_deRham_Laplacian}
\end{equation}
(e.g., see~\cite[equation~(5.4)]{Rosenberg-97}).
From \eqref{eq:scaling_deRham_Laplacian}, we see
\begin{align*}
 \kappa_{\deRham} (M, E , g_{\rho} ) (s)
 &{}=
 \sum_{k=0}^{2n+1} (-1)^k k \dim H_{\deRham}^k (M, E)
 \\
 &\quad{}+ {\rm e}^{2 \rho s }
 \sum_{k=0}^{2n+1} (-1)^k k \big(
 \zeta \big(\Delta_{\deRham , g}^k \big) (s )
 -
 \dim H_{\deRham}^k (M, E)
 \big)
 .
\end{align*}
To derivate the above equation,
\begin{align*}
 \kappa_{\deRham} (M, E , g_{\rho} ) ' (s)
 & =
 2 \rho {\rm e}^{2 \rho s }
 \sum_{k=0}^{2n+1} (-1)^k k \big(
 \zeta \big(\Delta_{\deRham , g}^k \big) (s )
 -
 \dim H_{\deRham}^k (M, E)
 \big)
 \\
 & \quad
 + {\rm e}^{2 \rho s }
 \sum_{k=0}^{2n+1} (-1)^k k \big(
 \zeta \big(\Delta_{\deRham , g}^k \big) ' (s )
 \big)
 .
\end{align*}
To substitute $s = 0$,
since
\[
\zeta \big(\Delta_{\deRham , g}^k \big) (0 ) = 0
\]
on manifolds with dimension $2n+1$,
the Ray--Singer torsion is given by
\begin{equation}
 \log T_{\deRham} (M, E , g_{\rho})
 =
 \log T_{\deRham} (M, E , g)
 - \rho \sum_{k = 0 }^{2n+1} (-1)^k k \dim H_{\deRham} (M, E)
 .
 \label{eq:scaling_analytic_torsion_2}
\end{equation}
To substitute $M = K$, $E$ is the trivial bundle, $\rho = \log 2 $ and $g = g_{\mathrm{std}}$,
we obtain
\begin{equation}
 T_{\deRham} (K, \underline{\mathbb{C}} , 4 g_{\mathrm{std}})
 = 2^{2n+1} T_{\deRham} (K, \underline{\mathbb{C}} , g_{\mathrm{std}}) .
 \label{eq:scaling_analytic_torsion}
\end{equation}
By \eqref{eq:pre_Proposition_trivial_bundle} and \eqref{eq:scaling_analytic_torsion},
we conclude Proposition~\ref{prop:Ray--Singer_torsion_of_the_trivial_bundle_on_lens_spaces}.

\section{Proof of Corollary~\ref{cor:Rumin-and-Ray--Singer-analytic-torsion-on-lens-spaces}}
\label{sec:Ray--Singer_torsion_and_contact_analytic_torsion}

Since $ \alpha (\gamma) \in \mathrm{U}(r)$
is diagonalizable by a unitary matrix,
we have
\[
 E_{\alpha}
 =
 E_{u_1}
 \oplus
 \cdots
 \oplus
 E_{u_r}
 .
\]
We recall the Ray--Singer torsion on lens spaces,
\begin{Proposition}[{\cite{Ray-70}}]
 \label{prop:Ray-70}
 For $u$ $(= 1 ,\ldots , \mu -1)$,
 \[
 T_{\deRham} ( K , E_{u}, 4 g_{\mathrm{std}} )
 =
 \prod_{j=1}^{n+1} \big| {\rm e}^{ 2 \pi \ii u {\tau}_j / \mu } -1 \big|
 .
 \]
\end{Proposition}
From Proposition~\ref{prop:Ray-70},
Theorem~\ref{theo:compute-precisely-kappa-on-the-lens-spaces},
and Proposition~\ref{prop:Ray--Singer_torsion_of_the_trivial_bundle_on_lens_spaces},
we conclude
\begin{align*}
 T_{\deRham} ( K , E_{\alpha}, g_{\theta , J} )
 & =
 \prod_{j=1}^r T_{\deRham} ( K , E_{u_j}, g_{\theta , J} )
 =
 \prod_{j=1}^r n!^{-\dim H^0 (K,E_{u_j})} T_{\RN} ( K , E_{u_j}, g_{\theta , J} )
 \\
 & =
 n!^{-\dim H^0 (K,E_{\alpha})} T_{\RN} ( K , E_{\alpha}, g_{\theta , J} )
 .
\end{align*}

\appendix

\section{Alternative shorter derivation of Proposition~\ref{prop:Ray--Singer_torsion_of_the_trivial_bundle_on_lens_spaces}}
The decomposition
\[
 \Omega^{\bullet} \big(S^{2n+1} , \underline{\mathbb{C}}\big)
 =
 \bigoplus_{u=0}^{\mu -1}
 \Omega^{\bullet} ( K , E_u )
\]
gives
\[
 \kappa_{\deRham} \big(S^{2n+1 } , \underline{\mathbb{C}} , 4 g_{\mathrm{std}}\big)
 =
 \sum_{u=0}^{\mu -1}
 \kappa_{\deRham} (K , E_u , 4 g_{\mathrm{std}})
 .
\]
In particular,
\[
 T_{\deRham} \big(S^{2n+1 } , \underline{\mathbb{C}} , 4 g_{\mathrm{std}}\big)
 =
 T_{\deRham} (K , \underline{\mathbb{C}} , 4 g_{\mathrm{std}})
 \prod_{u=1}^{\mu -1}
 T_{\deRham} (K , E_u , 4 g_{\mathrm{std}})
 .
\]
Combining this with \eqref{eq:scaling_analytic_torsion_2}, \eqref{eq:scaling_analytic_torsion}, Propositions~\ref{prop:Weng-You} and~\ref{prop:Ray-70} and~\cite{Melo-Spreafico-09}, we obtain
\[
 \frac{(4 \pi )^{n+1} }{ n!}
 =
 T_{\deRham} (K , \underline{\mathbb{C}} , 4 g_{\mathrm{std}})
 \prod_{u=1}^{\mu -1}
 \prod_{j=1}^{n+1}
 \big|
 {\rm e}^{2 \pi \sqrt{-1} u \tau_j / \mu } -1
 \big|
 =
 T_{\deRham} (K , \underline{\mathbb{C}} , 4 g_{\mathrm{std}})
 \mu^{n+1},
\]
whence Proposition~\ref{prop:Ray--Singer_torsion_of_the_trivial_bundle_on_lens_spaces}.

\subsection*{Acknowledgement}
The author is grateful to his supervisor Professor Kengo Hirachi for introducing this subject and for helpful comments.
This work was supported by the program for Leading Graduate Schools, MEXT, Japan.
The author also thanks the referees for their valuable comments.


\pdfbookmark[1]{References}{ref}

\LastPageEnding

\end{document}